\documentclass [leqno] {amsart}
\usepackage [cp1250] {inputenc}
\usepackage[T1]{fontenc}
\usepackage {amsfonts}
\usepackage {amsmath}
\usepackage {amsfonts}
\usepackage {amsthm}
\usepackage{amssymb,latexsym}
\usepackage{enumerate}
\numberwithin {equation} {section}
\DeclareMathOperator{\Jac}{Jac}
\DeclareMathOperator{\Hess}{Hess}
\DeclareMathOperator{\Tr}{Tr}
\newtheorem {thm}{Theorem}[section]
\newtheorem {lem}[thm]{Lemma}

\newcommand{\der}[2]{\frac{\partial #1}{\partial #2}}
\newcommand{\ld}{\ldots}
\newcommand{\druga}[3]{\frac{\partial^2 #1}{\partial #2 \partial #3}}
\newcommand{\Cn}{\mathbb{C}^n}
\begin {document}
\title{Irreducibility of the symmetric Yagzhev's maps}
\author{Sławomir Bakalarski}
\date{}
\maketitle
\begin {abstract}
Let $F:\Cn \rightarrow \Cn$ be a polynomial mapping in Yagzhev's form,i.e. $$F(x_1,\ld,x_n)=(x_1+H_1(x_1,\ld,x_n),\ld,x_n+H_n(x_1,\ld,x_n)),$$ where $H_i$ are homogenous polynomials of degree 3. In this paper we show that if $\Jac(F) \in \mathbb{C}^*$ and the Jacobian matrix of $F$ is symmetric, then all the polynomials $x_i+H_i(x_1,\ld,x_n)$ are irreducible as elements of the ring $\mathbb{C}[x_1,\ld,x_n]$.
\end {abstract}
\section {Introduction}
Let $F:\Cn \rightarrow \Cn$ be a polynomial mapping. The Jacobian Conjecture says that if $\Jac(F) \in \mathbb{C}^*$, then $F$ is a polynomial automorphism. In the study of this problem many authors have showed that it is sufficient to check it only for some special class of mappings. One of this special class form mappings in so called Yagzhev's form (\cite{BCW},\cite{D},\cite{Y}). The mapping $F:\Cn \rightarrow \Cn$ is as above if there exist homogenous polynomials $H_1,\ld,H_n$ of degree 3 such that $$F(x_1,\ld,x_n)=(x_1+H_1(x_1,\ld,x_n),\ld,x_n+H_n(x_1,\ld,x_n)).$$We treat zero polynomial as a homogenous polynomial of any degree. In \cite{BE},\cite{M} it is shown that one can go further and assume additionally that the Jacobian matrix of $F$ is also symmetric.

One of the questions related to the Jacobian Conjecture is a question whether the coordinates of the jacobian mapping (i.e. one with constant jacobian) are irreducible polynomials.
Of course if the Jacobian Conjecture was true, then the answer to the above question would be in affirmative. In the paper we prove that if the mapping $F:\Cn \rightarrow \Cn$ is in Yagzhev's form and has symmetric Jacobian matrix, then all its coordinates are irreducible as elements of the ring $\mathbb{C}[x_1,\ld,x_n]$.

\section {Main Result}
Let $F=(F_1,\ld,F_n):\Cn \rightarrow \Cn$ be a polynomial mapping in the Yagzhev's form. This means that there exist homogenous polynomials $H_1,\ld,H_n$ of degree 3 such that $$F(x_1,\ld,x_n)=(x_1+H_1(x_1,\ld,x_n),\ld,x_n+H_n(x_1,\ld,x_n)).$$ Let us define $f_k=x_k+H_k(x_1,\ld,x_n)$ for $k=1,\ld,n$. Now we shall prove a technical lemma which shows all the possible factorizations of the polynomial $f_k=x_k+H_k(x_1,\ld,x_n)$ if $H_k$ is not zero. 
\begin {lem}
With the above notations we have that if $H_k \neq 0$, then the polynomial $f_k$ is reducible if and only if $x_k$ divides $H_k$.
\end {lem}

\begin {proof}
One implication is obvious. If $f_k$ is reducible, then also its homogenization is reducible. Hence there exist homogenous polynomials $P_i$ of degree $i$, where $i=1,2$ such that $$(*) \; x_kx_{n+1}^2+H_k(x_1,\ld,x_n)=P_1P_2.$$ Write $P_1=(\alpha x_{n+1}+L_1(x_1,\ld,x_n))$ and $P_2=(\beta x_{n+1}^2+L_2(x_1,\ld,x_n)x_{n+1}+Q(x_1,\ld,x_n))$, where $L_1,L_2$ are linear polynomials and $Q$ is homogenous polynomial of degree 2. Comparing coefficients we obtain following equalities
\begin {enumerate}[(a)]
 \item $\alpha \beta=0$.
\item $\alpha L_2+\beta L_1=x_k$.
\item $\alpha Q+ L_1L_2=0$.
\item $L_1Q=H_k$
\end {enumerate}
If $\alpha=0$, then from (b) we have $L_1=\frac{1}{\beta}{x_k}$. According to $(*)$ we obtain that $x_k|H_k$. Now, if $\alpha \neq 0$ we get  $\beta=0$ and from (b) we obtain that  $L_2=\frac{1}{\alpha}x_k$. Substituting this to $(c)$ we get $$\alpha Q+\frac{1}{\alpha}L_1x_k=0.$$ This however shows that $x_k|Q$ and from (d) we have that $x_k|H_k$. This completes the proof.

\end {proof}
In the rest of the paper we need a closer characterization of mappings in the Yagzhev's form with the symmetric matrix $JF$. In \cite{BE},\cite{M} it is shown that the mapping  $F=(F_1,\ld,F_n):\Cn \rightarrow \Cn$ is as above if and only if there exists a homogenous polynomial $P \in \mathbb{C}[x_1,\ld,x_n]$ of degree 4 such that  $F_i=x_i+\der{P}{x_i}$. 

\bigskip
Now we have got all the facts needed to prove the main theorem.

\begin {thm}
Let $P$ be a homogenous polynomial of degree 4 and let \\ $F=(F_1,\ldots,F_n):\Cn \rightarrow \Cn$ be a polynomial mapping such that  $F_i=x_i+\frac{\partial P}{\partial x_i}$. If  $\Jac F=const \neq 0$, then the polynomials  $F_1,\ldots,F_n$ are irreducible.
\end {thm}

\begin {proof}
First notice that if $F$ is as above, then $\Jac(F) \in \mathbb{C}^*$ implies that $\Jac(F)=1$. 

We show that the assumption that some of the polynomials $F_i$ are reducible leads to a contradiction.

\textbf{1}. Without loss of the generality we can assume that the reducible polynomial is $F_n$. Therefore 
$\der{P}{x_n}$ is not zero and it is divisible by $x_n$; i.e.  $\der{P}{x_n}=x_n r(x_1,\ld,x_n)$ for some non-zero polynomial  $r$. Hence there exist polynomials $g,h$ such that
$$P=x_n^2g(x_1,\ld,x_n)+h(x_1,\ld,x_{n-1}),$$ where $g$ is 
a non-zereo, homogenous polynomial of degree 2 and $h$ is a homogenous polynomial of degree 4 in variables $x_1,\ld,x_{n-1}$.

Taking this into account, we obtain that the form of $F$ is as follows:
 $$F_i(x_1,\ldots,x_n)=x_i+\frac{\partial P}{\partial x_i}=x_i+x_n^2(a_{i1}x_1+\ldots+a_{in}x_n)+\der{h}{x_i},$$ for some $a_{i1},\ldots,a_{in}\in \mathbb{C}$, $i=1,\ld,n-1$ and
$$F_n(x_1,\ldots,x_n)=x_n+\frac{\partial P}{\partial x_n}=x_n+x_nQ(x_1,\ldots,x_n),$$ where $Q$ is a non-zero homogenous polynomial of degree 2. Write
$$Q(x_1,\ldots,x_n)=\sum^{n}_{i=1}c_{ii}x_i^2+\sum_{i<j}c_{ij}x_ix_j$$ and put $A=[a_{ij}]_{i,j=1,\ldots,n-1}$.

Since the matrix  $JF$ is symmetric, then in particular for every $k,l \in \{1,\ld,n-1\}$ we have equality $$\der{F_k}{x_l}=\der{F_l}{x_k}. $$ 
By comparison of coefficients at $x_n^2$, we obtain that the matrix $A$ is symmetric.

\textbf{2}. Using the fact that the matrix $JF$ is also symmetric we obtain that $$\der{F_1}{x_n}=\der{F_n}{x_1},$$ thus
$$2x_n(a_{11}x_1+\ldots+a_{1,n-1}x_{n-1})+3a_{1n}x_n^2=2c_{11}x_1x_n+\sum_{1<j\leq n}c_{1j}{x_j}x_n;$$ therefore
$$c_{11}=a_{11},\; c_{1j}=2a_{1j}, \; j=2,\ld,n-1, \; c_{1n}=3a_{1n}.$$
Analogously, from the fact that  $\der{F_l}{x_n}=\der{F_n}{x_l}$ for $l=2,\ld,n-1$ we get
$$c_{ll}=a_{ll},\; c_{lj}=2a_{lj}\text{ for } j=l+1,\ld,n-1,\;\;  c_{ln}=3a_{ln}.$$

\textbf{3}. Let us write $Q(x_1,\ld,x_n)=R(x_1,\ld,x_{n-1})+x_n\tilde{Q}(x_1,\ld,x_n)$, where $R(x_1,\ld,x_{n-1})$ does not depend on the variable  $x_n$. Hence we have the equality  $$\frac{\partial}{\partial x_n}(x_n Q)=R(x_1,\ld,x_{n-1})+x_n q(x_1,\ld,x_n)$$ for some polynomial  $q$.  \\ Therefore
$JF$ is equal to 

\smallskip
\begin{tiny}
$
\begin {pmatrix}
 1+a_{11}x_n^2+\frac{\partial^2 h}{\partial x_1^2} & a_{12}x_n^2+\druga{h}{x_2}{x_1} & \hdots & a_{1,n-1}x_n^2+\druga{h}{x_{n-1}}{x_1} & 3a_{1n}x_n^2+2x_n(a_{11}x_1+\ld+a_{1,n-1}x_{n-1}) \\  
a_{21}x_n^2+\druga{h}{x_1}{x_2} & 1+a_{22}x_n^2+\frac{\partial^2 h}{\partial x_2^2}& \hdots  & a_{2,n-1}x_n^2+\druga{h}{x_{n-1}}{x_2} & 3a_{2n}x_n^2+2x_n(a_{21}x_1+\ld+a_{2,n-1}x_{n-1}) \\
\vdots             & \vdots                  & \vdots  & \vdots & \vdots \\
a_{n-1,1}x_n^2+\druga{h}{x_1}{x_{n-1}} & a_{n-1,2}x_n^2+\druga{h}{x_2}{x_{n-1}} & \hdots   & 1+a_{n-1,n-1}x_n^2+\frac{\partial^2 h}{\partial x_{n-1}^2}& 3a_{n-1,n}x_n^2+2x_n(a_{n-1,1}x_1+\ld+a_{n-1,n-1}x_{n-1}) \\  
x_n \der{Q}{x_1} & x_n \der{Q}{x_2} & \hdots & x_n \der{Q}{x_{n-1}} & 1+R(x_1,\ld,x_{n-1})+x_nq(x_1,\ld,x_n)
\end {pmatrix}
$
\end{tiny} 
We substitute $x_n=0$ to the jacobian and obtain the following equality
\begin {equation*}
 1=\det JF|_{x_n=0}=\det 
\begin {pmatrix}
* & 0\\
0         & 1+R(x_1,\ld,x_{n-1})
\end {pmatrix}.
\end {equation*}
By the Laplace formula the above gives us that $\det(*)$ as well as $1+R(x_1,\ld,x_{n-1})$ are non-zero constants. If $1+R(x_1,\ld,x_{n-1}) \neq 1$, then  
 $R(x_1,\ld,x_{n-1})$ would be a non-zero constant contradicting the fact that $R(0,\ld,0)=0$. Therefore we have  $\det(*)=1$ and $1+R(x_1,\ld,x_{n-1})=1$, so $R(x_1,\ld,x_{n-1})=0$ and  $c_{ij}=0$ for all $i,j=1,\ld,n-1$. According to step \textbf{2} and the symmetry of the matrix $A$, we have that $A=0$.

\textbf{4}. Due to above considerations we obtain that the mapping $F$ is of the following form

\begin {equation*}
 F(x_1,\ld,x_n)=\begin {pmatrix}
 x_1+a_{1n}x_n^3+\der{h}{x_1} \\
x_2+a_{2n}x_n^3+\der{h}{x_2} \\
\hdots \\
x_{n-1}+a_{n-1,n}x_n^3+\der{h}{x_{n-1}} \\
x_n+x_n(3a_{1n}x_1x_n+\ld+3a_{n-1,n}x_{n-1}x_n+c_{nn}x_n^2)
\end {pmatrix}.
\end {equation*}
Since $F_i=x_i+\der{P}{x_i}$ we have 
$$P=h+a_{1n}x_1x_n^3+a_{12}x_2x_n^3+\ld+a_{n-1,n}x_{n-1}x_n^3+\frac{1}{4}c_{nn}x_n^4.$$

\textbf{5}.  The matrix $\Hess(P)$ of the hessian of the polynomial $P$ is nilpotent because the jacobian of $F$ is non-zero constant (\cite{D}), in particular its trace is equal to zero. By the above formula for $P$ we get

$$ 0=\Tr(\Hess(P))=\sum^{n-1}_{i=1}\frac{\partial^2h}{\partial x_i^2}+\frac{\partial^2P}{\partial x_n^2}=$$
$$ =w(x_1,...,x_{n-1})+6a_{1n}x_1x_n+6a_{2n}x_2x_n+\ld+6a_{n-1,n}x_{n-1,n}x_n+
3c_{nn}x_n^2.$$

Therefore
$$6a_{1n}x_1x_n+6a_{2n}x_2x_n+\ld+6a_{n-1,n}x_{n-1,n}x_n+3c_{nn}x_n^2=0,$$
thus

$$a_{1n}=a_{2n}=\ld=a_{n-1,n}=c_{nn}=0.$$ In this way we obtain that $Q=0$-this contradicts that  $Q$ is a non-zero homogenous polynomial of degree 2. This proves that $F_n$ is irreducible and completes the proof.

\end {proof}

\bigskip
\curraddr{
Jagiellonian University \\
Institute of Mathematics \\
Reymonta 4 \\
30-059 Kraków \\ 
Poland\\
E-mail: Slawomir.Bakalarski@im.uj.edu.pl
}


\begin {thebibliography}{[****]}
\bibitem [BCW]{BCW} H. Bass, E. Connell, D. Wright, {\em The Jacobian Conjecture: reduction of degree and formal expansion of the inverse}, Bull. Amer. Math. Soc. \textbf{7} (1982), 287-330.
\bibitem [BE]{BE} M. de Bondt, A. van den Essen, {\em A reduction of the Jacobian Conjecture to the symmetric case}, Proc. Amer. Math. Soc. , \textbf{133} (2005), no. 8, 2201-2205.
 \bibitem [D]{D} L. M. Drużkowski, {\em An effective approach to Keller's Jacobian Conjecture}, Math. Ann. \textbf{264} (1983), 303-313.
\bibitem [M] {M} G. Meng, {\em Legendre Transform, Hessian Conjecture and Tree Formula}, Appl. Math. Lett. 19 (2006), 503-510;   arXiv:math-ph/0308035v2.
\bibitem [Y]{Y} A. Yagzhev, {\em On Keller's problem}, Siberian Math. J. 21 (1980), 747-754.





\end {thebibliography}
\end {document}